\documentclass[11pt]{amsart}
\usepackage[T1]{fontenc}
\usepackage[utf8]{inputenc}

\usepackage{amsmath,amssymb,amsthm,amsfonts,mathtools,bm}

\usepackage{lmodern}
\usepackage[final]{microtype}
\usepackage{babel}
\usepackage{nicefrac} 

\usepackage[margin=1.2in]{geometry}

\usepackage{multicol}
\usepackage{array}
\usepackage{booktabs}
\usepackage[shortlabels]{enumitem}

\usepackage{tikz}
\usetikzlibrary{cd}
\usetikzlibrary{shapes}
\usetikzlibrary{shapes.multipart}

\usepackage{thmtools}
\usepackage{thm-restate}
\usepackage[hyphens]{url}
\usepackage[colorlinks]{hyperref}
\hypersetup{
    colorlinks,
    linkcolor={red!50!black},
    citecolor={blue!50!black},
    urlcolor={blue!80!black}
}

\theoremstyle{plain}
\newtheorem{theorem}{Theorem}[section]

\newtheorem*{claim*}{Claim}

\newtheorem*{proposition*}{Proposition}
\newtheorem{fact}[theorem]{Fact}
\newtheorem*{fact*}{Fact}

\newtheorem*{conjecture*}{Conjecture}

\newtheorem{lemma}[theorem]{Lemma}
\newtheorem*{lemma*}{Lemma}

\newtheorem*{question*}{Question}
\theoremstyle{definition}\newtheorem{remark}[theorem]{Remark}
\theoremstyle{definition}\newtheorem*{remark*}{Remark}
\theoremstyle{definition}\newtheorem{definition}[theorem]{Definition}
\theoremstyle{definition}\newtheorem*{definition*}{Definition}
\theoremstyle{definition}
\theoremstyle{definition}
\theoremstyle{definition}
\theoremstyle{definition}\newtheorem*{example*}{Example}

\newcommand{\term}{\textbf} 
\newcommand{\defeq}{\vcentcolon=} 
\newcommand{\concat}{^\smallfrown} 
\newcommand{\forces}{\Vdash} 
\newcommand{\0}{\varnothing} 
\renewcommand{\phi}{\varphi} 
\renewcommand{\epsilon}{\varepsilon} 

\def\N{\mathbb{N}} 
\def\Cantor{2^\N} 

\def\LO{\mathsf{LO}}
\DeclareMathOperator{\dist}{dist}
\DeclareMathOperator{\diam}{diam}
\def\Strings{2^{<\N}}
\def\BSigma{\mathbf{\Sigma}}

\title[Borel graphability and diameter]{A Borel graphable equivalence relation with no Borel graphing of diameter two}
\author{Patrick Lutz}
\address{Department of Mathematics, University of Michigan}
\email{pglutz@umich.edu}

\begin{document}

\begin{abstract}
We answer a question of Arant, Kechris and Lutz by showing that there is a Borel graphable equivalence relation with no Borel graphing of diameter less than $3$. More specifically, we prove that there is an equivalence relation with a Borel graphing of diameter at most $4$ but no Borel graphing of diameter less than $3$. Our proof relies on a technical lemma about computability-theoretic genericity, which may have other applications.
\end{abstract}

\maketitle

\vspace{-10pt}\section{Introduction}

Suppose that $E$ is an equivalence relation on a Polish space $X$. A \term{Borel graphing} of $E$  is a Borel graph $G$ on $X$ such that for all $x, y \in X$,
\[
  x\, E\, y \iff \text{ there is a path in $G$ between $x$ and $y$}.
\]
In other words, the connected components of $G$ are exactly the equivalence classes of $E$. We say that $E$ is \term{Borel graphable} if it has a Borel graphing.

The \term{diameter} of a graphing $G$ of $E$ is defined as the supremum of the distances in $G$ over all pairs of $E$-equivalent\footnote{Note that this slightly differs from the usual definition of the diameter of a graph, in which we would take the supremum over all pairs of elements, not just pairs that are $E$-equivalent. Our definition is equivalent to taking the supremum of the diameters (in the standard sense) of the connected components of the graph.} elements of $X$. More precisely, for any pair $x, y \in X$ such that $x\, E\, y$, we use $\dist_G(x, y)$ to denote the length of the shortest path\footnote{Note that we use the convention that if $x = y$ then there is always a path of length $0$ between $x$ and $y$.} in $G$ between $x$ and $y$. We then define the diameter of $G$ to be
\[
  \diam(G) \defeq \sup\{\dist_G(x, y) \mid x\, E\, y\}
\]
where if $\{\dist_G(x, y) \mid x\, E\, y\}$ is unbounded in $\N$ then $\diam(G) = \infty$. 

The notion of a Borel graphable equivalence relation was introduced by Arant~\cite{arant2019effective} and further investigated by Arant, Kechris and Lutz~\cite{arant2024borel}. In the latter paper, it was observed that when an equivalence relation is shown to be Borel graphable, it is nearly always done by exhibiting a Borel graphing of diameter at most $2$. In light of this, it was asked in that paper whether there is a Borel graphable equivalence relation for which this cannot be done. More precisely, is there a Borel graphable equivalence relation $E$ such that every Borel graphing of $E$ has diameter at least $3$? The goal of this paper is to answer this question in the affirmative.

\begin{theorem}
\label{thm:main}
There is an equivalence relation on a Polish space which is Borel graphable but for which every Borel graphing has diameter at least $3$.
\end{theorem}

In order to prove this theorem, we will define a specific equivalence relation and show that it has a Borel graphing of diameter at most $4$, but no Borel graphing of diameter less than $3$ (incidentally, we do not know whether it has a Borel graphing of diameter exactly $3$, though we suspect it does not). In fact, this equivalence relation is actually one that was defined in the paper by Arant, Kechris and Lutz mentioned above, though in order to keep this paper relatively self-contained we review the definition in Section~\ref{sec:theorem}. In that paper, it was also shown that this equivalence relation has a Borel graphing of diameter at most $4$; our contribution lies strictly in showing that it has no Borel graphing of diameter less than $3$.

The definition of this equivalence relation uses the notion of computability-theoretic genericity, which we will review in Section~\ref{sec:prelim}. A key tool in our proof is the following technical lemma about this sort of genericity, which we believe could have other applications.

\begin{lemma}
\label{lemma:key}
For any $r \in \Cantor$, there are $x, y \in \Cantor$ such that $x$ and $y$ are $r'$-generic and for any $z\in \Cantor$ which is $r'$-generic, either $x$ and $z$ are mutually $r$-generic or $y$ and $z$ are mutually $r$-generic.
\end{lemma}

\begin{remark}
Our proof of this lemma can be easily modified to prove a slightly stronger statement: for all $r, s \in \Cantor$, there are $x, y \in \Cantor$ such that $x$ and $y$ are $s$-generic and for any $z \in \Cantor$ which is $r'$-generic, either $x$ and $z$ are mutually $r$-generic or $y$ and $z$ are mutually $r$-generic. In other words, $x$ and $y$ can be required to be as generic as we like (i.e.\ much more than just $r'$-generic) without affecting the conclusion of the lemma.
\end{remark}

In Section~\ref{sec:theorem}, we will see how to use this lemma to prove Theorem~\ref{thm:main}, before proving the lemma itself in Section~\ref{sec:lemma}. Before moving on to the technical details of this paper, however, we would like to note that there are a few questions about diameters of Borel graphings which we do not know the answer to and for which it is not clear if the techniques of this paper are helpful.

First, for each $n \in \N$, one can ask if there is an equivalence relation with a Borel graphing of diameter $n + 1$, but no Borel graphing of diameter $n$. Note that for $n = 1$, this is equivalent to asking if there is an equivalence relation which has a Borel graphing of diameter $2$ but which is not Borel and that many such equivalence relations were shown to exist in~\cite{arant2019effective} and~\cite{arant2024borel}. As we mentioned above, the equivalence relation we use in the proof of Theorem~\ref{thm:main} has a Borel graphing of diameter $4$, but no Borel graphing of diameter of diameter $2$. Hence this equivalence relation provides a positive answer to the question for either $n = 2$ or $n = 3$, though we don't know which (but we suspect for $n = 3$). For all other values of $n$, however, we do not know the answer to this question and the ideas used in our proof of Theorem~\ref{thm:main} do not seem immediately helpful.

A second, and closely related, question is whether there is a Borel graphable equivalence relation with no Borel graphing of finite diameter. Note that if we had a positive answer to the first question for infinitely many values of $n$ then we would also have a positive answer to this question. In particular, if we had a sequence of Borel graphable equivalence relations $\{E_n\}_n$ where each $E_n$ has no Borel graphing of diameter less than $n$ then we could take their disjoint union to get a Borel graphable equivalence relation with no Borel graphing of finite diameter. However, for the converse direction the implication is not clear. And once again, the ideas used in the proof of Theorem~\ref{thm:main} do not seem immediately helpful in answering this question.

\section{Background on computability-theoretic genericity}
\label{sec:prelim}

We will now review the standard computability-theoretic definition of genericity. In what follows, we will assume familiarity with effective descriptive set theory, forcing, and the hyperarithmetic hierarchy; for an introduction to these topics, see the book \emph{Recursive Aspects of Descriptive Set Theory} by Mansfield and Weitkamp~\cite{mansfield1985recursive}. For a more thorough introduction to computability-theoretic genericity, see Section 2.24 of \emph{Algorithmic Randomness and Complexity} by Downey and Hirschfeldt~\cite{downey2010algorithmic}.

\subsection{Genericity and mutual genericity}

Roughly speaking, given $r \in \Cantor$, $x \in \Cantor$ is $r$-generic if it is generic for Cohen forcing when we consider only those dense sets which are computable from $r$. However, the precise definition is a bit more complicated than this.\footnote{Actually, the informal gloss on computability-theoretic genericity that we have given here is usually referred to in computability theory as ``weak genericity.'' In particular, $x \in \Cantor$ is \term{weakly $r$-generic} if it meets every dense set for Cohen forcing which is computable from $r$. For the purposes of this paper, the difference between genericity and weak genericity is not especially important; the only reason we have chosen to use genericity rather than weak genericity is that it is the more standard notion in computability theory.}

\begin{definition}
Suppose $A \subseteq \Strings$ is a set of strings. A real $x \in \Cantor$ \term{meets} $A$ if there is some $n \in \N$ such that $x\restriction n \in A$ and \term{avoids} $A$ if there is some $n \in \N$ such that for all $\sigma \geq x\restriction n$, $\sigma \notin A$.
\end{definition}

\begin{definition}[Jockusch and Posner~\cite{jockusch1980degrees}]
Given some fixed $r \in \Cantor$, $x \in \Cantor$ is \term{$r$-generic} if for every set $A \subseteq \Strings$ which is c.e.\ relative to $r$, $x$ either meets or avoids $A$.
\end{definition}

It is not hard to see that for any $r$, there are continuum-many $r$-generics. It is also not hard to see that as $r$ increases in Turing degree, the requirement of $r$-genericity becomes more stringent. In particular, we have the following facts, which we will make use of below.

\begin{fact}
For every $r, s \in \Cantor$, if $r \leq_T s$ then every element of $\Cantor$ which is $s$-generic is also $r$-generic.
\end{fact}

\begin{fact}
\label{fact:generic_not_generic}
For every $r \in \Cantor$, there is some $x \in \Cantor$ such that $x$ is $r$-generic, but not $r'$-generic.
\end{fact}

There is also a computability-theoretic notion of mutual genericity: $x$ and $y$ are \term{mutually $r$-generic} if $x\oplus y$ is $r$-generic (where $x\oplus y$ denotes the standard way of encoding two elements of $\Cantor$ as a single element of $\Cantor$). It is fairly easy to see that this can be redefined in terms of meeting or avoiding c.e.\ sets of pairs of strings.

\begin{definition}
Suppose that $A \subseteq \Strings\times \Strings$ is a set of pairs of strings. A pair of reals $(x, y) \in \Cantor\times \Cantor$ \term{meets} $A$ if there are some $n, m \in \N$ such that $(x\restriction n, y\restriction m) \in A$ and \term{avoids} $A$ if there are some $n,m \in \N$ such that for all $\sigma \geq x\restriction n$ and $\tau \geq y\restriction m$, $(\sigma, \tau) \notin A$.
\end{definition}

\begin{fact}
For any $r \in \Cantor$ and $x, y \in \Cantor$, $x$ and $y$ are mutually $r$-generic if and only if for every set of pairs of strings $A \subseteq \Strings\times\Strings$ which is c.e.\ relative to $r$, the pair $(x, y)$ either meets or avoids $A$.
\end{fact}

Another characterization of mutual genericity, which we will use below, was proved by Liang Yu~\cite{yu2006lowness}. It is sometimes referred to as the analogue of van Lambalgen's theorem for genericity.\footnote{van Lambalgen's theorem is an important theorem in the theory of algorithmic randomness which inspired Yu's theorem.}

\begin{theorem}[Yu]
\label{thm:van_lambalgen}
For any $r \in \Cantor$, $x$ and $y$ are mutually $r$-generic if and only if $x$ is $r$-generic and $y$ is $(r\oplus x)$-generic.
\end{theorem}

\subsection{Genericity and forcing}
\label{sec:generic_forcing}

As indicated by both the name and some of our comments above, the computability-theoretic notion of genericity is closely connected to Cohen forcing. There are two points about this connection that will be useful later.

First, the definition of genericity can be rephrased in terms of meeting certain dense sets for Cohen forcing. Recall that a set of strings $D \subseteq \Strings$ is \term{dense} if for all $\sigma \in \Strings$, there is some $\tau \geq \sigma$ such that $\tau \in D$. For any $r \in \Cantor$, $r$-genericity can be characterized in terms of meeting certain dense sets depending on $r$. In particular, given any set of strings $A \subseteq \Strings$, let $D_A$ be the set of strings defined by
\[
  D_A = \{\sigma \mid \sigma \in A \text{ or } \forall \tau \geq \sigma\, (\tau \notin A)\}.
\]
It is easy to check that $D_A$ is dense and that $x$ is $r$-generic if and only if $x$ meets each dense set of the form $D_A$, where $A$ is c.e.\ relative to $r$.

A simple, but useful, point is that if $D$ is a set of strings which is computable from $r$ then $D$ can be written in the form $D_A$ for some $A$ which is c.e.\ relative to $r$ (in fact, we can just take $A = D$). Thus if $x$ is $r$-generic then it meets every dense set which is computable from $r$.

Of course, nearly identical comments apply for mutual genericity. In particular, a set of pairs of strings $D \subseteq \Strings\times\Strings$ is \term{dense} if for every pair of strings $(\sigma, \tau) \in \Strings\times\Strings$, there is a pair of strings $(\sigma', \tau') \in D$ such that $\sigma' \geq \sigma$ and $\tau' \geq \tau$. To any set of pairs of strings $A \subseteq \Strings \times \Strings$ we can associate a dense set $D_A$ such that $x$ and $y$ are mutually $r$-generic if and only if the pair $(x, y)$ meets every dense set of the form $D_A$ for $A$ c.e.\ relative to $r$.

Second, if $x$ is $r$-generic then $\Sigma^0_1(r)$ and $\Pi^0_1(r)$ facts about $x$ must be forced by some initial segment of $x$ (in fact, this actually characterizes $r$-genericity and is the motivation for the slightly complicated definition). In particular, if $x$ is $r$-generic then for every $\Sigma^0_1$ formula $\phi$ with $r$ as a parameter, there is some initial segment $\sigma$ of $x$ such that either $\sigma \forces \phi$ or $\sigma \forces \lnot \phi$, where $\forces$ denotes the forcing relation for Cohen forcing. One consequence of this is that if $x$ is $r$-generic and is contained in some $\Sigma^0_1(r)$ set $A \subseteq \Cantor$ then there is some initial segment $\sigma$ of $x$ such that for any $r$-generic $y$ which extends $\sigma$, $y$ is contained in $A$ as well. Moreover, this fact can be generalized to all levels of the lightface Borel hierarchy as follows.

\begin{fact}
\label{fact:generic_forces}
Fix $r \in \Cantor$ and $\alpha < \omega_1^r$ and suppose that $A\subseteq \Cantor$ belongs to the lightface pointclass $\Sigma^0_{1 + \alpha}(r)$. If $x$ is $r^{(\alpha)}$-generic and contained in $A$ then there is some finite initial segment $\sigma$ of $x$ such that any $r^{(\alpha)}$-generic $y$ which extends $\sigma$ is also contained in $A$.
\end{fact}

Once again, nearly identical facts are true of mutual genericity. In particular, we have the following fact.

\begin{fact}
\label{fact:mutually_generic_forces}
Fix $r \in \Cantor$ and $\alpha < \omega_1^r$ and suppose that $A\subseteq \Cantor\times \Cantor$ belongs to the lightface pointclass $\Sigma^0_{1 + \alpha}(r)$. If $x$ and $y$ are mutually $r^{(\alpha)}$-generic and $(x, y) \in A$ then there are finite initial segments $\sigma$ of $x$ and $\tau$ of $y$ such that for any mutually $r^{(\alpha)}$-generic $u$ and $v$ with $u$ extending $\sigma$ and $v$ extending $\tau$, $(u, v) \in A$.
\end{fact}

\section{Proof of the main theorem}
\label{sec:theorem}

We will now show that there is an equivalence relation which is Borel graphable but has no Borel graphing of diameter less than $3$. We will begin by defining the equivalence relation, which, as we mentioned above, was first studied by Arant, Kechris and Lutz in~\cite{arant2024borel}.

Let $\LO$ denote the Polish space of linear orders with domain $\N$ (see Sections 16C and 27C of~\cite{kechris1995classical} for a more detailed definition). For a well-founded linear order $L \in \LO$, let $|L|$ denote the ordinal that $L$ is isomorphic to. Let $X$ be the Polish space $\LO\times\Cantor\times\Cantor$ and let $E$ be the equivalence relation on $X$ defined by setting $(L, r, x)$ and $(R, s, y)$ equivalent if $L = R$, $r = s$, and one of the following holds:
\begin{enumerate}
\item $L$ is ill-founded.
\item $L$ is well-founded and neither $x$ nor $y$ is $r^{(\alpha)}$-generic, where $\alpha = |L|$.
\item $L$ is well-founded and both $x$ and $y$ are $r^{(\alpha)}$-generic, where $\alpha = |L|$.
\end{enumerate}
In other words, for each $L \in \LO$ and $r \in \Cantor$, the set $\{(L, r, x) \mid x \in \Cantor\}$ consists of either one or two $E$-equivalence classes: one if $L$ is ill-founded and two if $L$ is well-founded, in which case one equivalence class consists of those tuples whose third coordinate is $r^{(\alpha)}$-generic and the other consists of those tuples whose third coordinate is not $r^{(\alpha)}$-generic (where $\alpha = |L|$). 

It was shown in~\cite{arant2024borel} that $E$ is analytic (Proposition~57) and Borel graphable by a Borel graphing with diameter at most $4$ (Proposition~60). Thus to prove Theorem~\ref{thm:main}, it is enough to prove that $E$ has no Borel graphing of diameter less than $3$.

Suppose for contradiction that $E$ does have a Borel graphing of diameter less than $3$ and let $G$ be such a graphing. Since $G$ is Borel, there is some $\alpha < \omega_1$ such that $G$ is in the boldface pointclass $\BSigma^0_{1 + \alpha}$ and hence there is some $r \in \Cantor$ such that $\alpha < \omega_1^r$ and $G$ is in the lightface pointclass $\Sigma^0_{1 + \alpha}(r)$. Fix such an $r$ and let $L$ be a presentation of $\alpha$ computable from $r$.

We will reason exclusively about elements of $X$ of the form $(L + 1, r, -)$, where $L + 1$ denotes some fixed presentation of $\alpha + 1$ which is computable from $r$. The overall strategy of the proof is as follows. First, we will show---in Lemma~\ref{lemma:small_key} below---that if $u$ and $v$ are mutually $r^{(\alpha)}$-generic then there cannot be an edge in $G$ between $(L + 1, r, u)$ and $(L + 1, r, v)$, essentially because $G$ is too simple to be able to tell if $u$ and $v$ are $r^{(\alpha + 1)}$-generic or not. Second, we will use Lemma~\ref{lemma:key} to show that $G$ must contain an edge of this form.

\begin{lemma}
\label{lemma:small_key}
Suppose that $u$ and $v$ are mutually $r^{(\alpha)}$-generic. Then there is no edge in $G$ between $(L + 1, r, u)$ and $(L + 1, r, v)$.
\end{lemma}

\begin{proof}
Note that since $G$ is $\Sigma^0_{1 + \alpha}(r)$, so is the relation $R \subseteq \Cantor\times \Cantor$ defined by
\[
  R(x, y) \iff \text{ there is an edge in $G$ between $(L + 1, x)$ and $(L + 1, r, y)$}.
\]
Now suppose for contradiction that there is an edge in $G$ between $(L + 1, r, u)$ and $(L + 1, r, v)$---i.e.\ that $R(u, v)$ holds. Since $u$ and $v$ are mutually $r^{(\alpha)}$-generic, by Fact~\ref{fact:mutually_generic_forces} there must be finite initial segments $\sigma$ of $u$ and $\tau$ of $v$ which force $R(u, v)$. In other words, for all $\tilde{u}$ extending $\sigma$ and $\tilde{v}$ extending $\tau$, if $\tilde{u}$ and $\tilde{v}$ are mutually $r^{(\alpha)}$-generic then $R(\tilde{u}, \tilde{v})$ holds.

To finish the proof, we can simply take some $\tilde{u}$ extending $\sigma$ and $\tilde{v}$ extending $\tau$ such that $\tilde{u}$ and $\tilde{v}$ are mutually $r^{(\alpha)}$-generic, $\tilde{u}$ is $r^{(\alpha + 1)}$-generic and $\tilde{v}$ is not $r^{(\alpha + 1)}$-generic.\footnote{For example, use Fact~\ref{fact:generic_not_generic} to find some $\tilde{v}$ extending $\tau$ which is $r^{(\alpha)}$-generic but not $r^{(\alpha + 1)}$-generic and then take any $\tilde{u}$ extending $\sigma$ which is $(\tilde{v}\oplus r^{(\alpha + 1)})$-generic. By Theorem~\ref{thm:van_lambalgen} (the genericity version of van Lambalgen's theorem), $\tilde{u}$ and $\tilde{v}$ are mutually $r^{(\alpha)}$-generic.} To see why this finishes the proof, note that by definition of $E$, $(L + 1, r, \tilde{u})$ and $(L + 1, r, \tilde{v})$ are not $E$-equivalent, but since $R(\tilde{u}, \tilde{v})$ holds, there is an edge between them in $G$, which contradicts the assumption that $G$ is a graphing of $E$.
\end{proof}

By Lemma~\ref{lemma:key}, we can find some $x, y \in \Cantor$ such that $x$ and $y$ are both $r^{(\alpha + 1)}$-generic and for any $z$ which is $r^{(\alpha + 1)}$-generic, either $x$ and $z$ are mutually $r^{(\alpha)}$-generic or $y$ and $z$ are mutually $r^{(\alpha)}$-generic. Note that $(L + 1, r, x)$ and $(L + 1, r, y)$ are $E$-equivalent and hence they must be connected in $G$ by a path of length either $1$ or $2$. We will now show that either of these two possibilities leads to a contradiction.

First, suppose that $(L + 1, r, x)$ and $(L + 1, r, y)$ are connected in $G$ by a path of length $2$. Let $(L + 1, r, z)$ be the middle point of this path---i.e.\ there are edges in $G$ between $(L + 1, r, x)$ and $(L + 1, r, z)$ and between $(L + 1, r, z)$ and $(L + 1, r, y)$. Since $(L + 1, r, z)$ is $E$-equivalent to $(L + 1, r, x)$, $z$ must be $r^{(\alpha + 1)}$-generic. Hence by our choice of $x$ and $y$, either $x$ and $z$ are mutually $r^{(\alpha)}$-generic or $y$ and $z$ are mutually $r^{(\alpha)}$-generic. Either way, we violate Lemma~\ref{lemma:small_key}.

Now suppose that $(L + 1, r, x)$ and $(L + 1, r, y)$ are connected in $G$ by a path of length $1$. In other words, there is an edge between them in $G$. We claim that $x$ and $y$ are mutually $r^{(\alpha)}$-generic, thus contradicting Lemma~\ref{lemma:small_key}. To see why, simply note that since $x$ is $r^{(\alpha + 1)}$-generic, then by taking $z = x$ in the statement of Lemma~\ref{lemma:key}, we have that either $x$ and $x$ are mutually $r^{(\alpha)}$-generic or that $x$ and $y$ are mutually $r^{(\alpha)}$-generic. As the former is impossible (an element of $\Cantor$ cannot be mutually generic with itself), the latter must hold.

\section{Proof of the key lemma}
\label{sec:lemma}




We will now prove Lemma~\ref{lemma:key}. Fix $r \in \Cantor$; our goal is to construct $r'$-generics $x$ and $y$ such that for any $r'$-generic $z$, either $x$ and $z$ are mutually $r$-generic or $y$ and $z$ are mutually $r$-generic.

We begin by fixing a bit of notation. Let $D_0, D_1, D_2,\ldots$ be an enumeration of the dense subsets of $2^{< \N}\times 2^{< \N}$ that a pair of reals must meet to be mutually $r$-generic (see Section~\ref{sec:generic_forcing}). Note that we can choose the sequence $D_0, D_1, D_2, \ldots$ to be uniformly computable\footnote{To see this, simply note that, in the notation of Section~\ref{sec:generic_forcing}, $r'$ can uniformly compute the dense set $D_A$ from an index for a set $A$ which is c.e.\ relative to $r$.} from $r'$.

We will now attempt to briefly explain some of the main ideas of the proof. A simple observation is that $x$ and $y$ must not be mutually $r'$-generic: if they were then we could take $z = x\oplus y$, which would be $r'$-generic by Theorem~\ref{thm:van_lambalgen} but not mutually $r$-generic with either $x$ or $y$. This suggests that our proof should follow other proofs in which two generic, but not mutually generic, reals are constructed (such as the proof of the well-known fact that for any $a \in \Cantor$, there are generic reals $b$ and $c$ such that $b\oplus c \geq_T a$). And this is indeed what we will do. Roughly speaking, $x$ and $y$ will trade off between meeting dense sets to ensure $r'$-genericity and providing opportunities for a potential $z$ to meet each $D_n$ with them. Slightly more precisely, at any given point in the construction, one of $x$ and $y$ will be trying to meet some dense set to ensure $r'$-genericity and the other will be trying to provide any potential $z$ with an opportunity to meet some $D_n$ together. The genericity of $z$ will then be used to ensure that $z$ takes one of these opportunities infinitely often.

The following lemma will help in ensuring the sort of coordination described above.

\begin{lemma}
\label{lemma:process}
There are $r'$-computable functions $g_1, g_2\colon \N \to 2^{< \N}$ such that for any $n$, and any strings $\gamma_1$ and $\gamma_2$ of length $n$, the pair $(\gamma_1\concat g_1(n), \gamma_2\concat g_2(n))$ meets all of the dense sets $D_0, D_1, \ldots, D_n$.
\end{lemma}

\begin{proof}
Fix $n$. We will define $g_1(n)$ and $g_2(n)$ via a sequence of $2^{2n}$ steps, one for each possible value of the pair $(\gamma_1, \gamma_2)$. More precisely, we will define two sequences of strings, $\0 = \sigma_0 \leq \sigma_1 \leq \sigma_2 \leq \ldots \leq \sigma_{2^{2n}}$ and $\0 = \tau_0 \leq \tau_1 \leq \tau_2 \leq \ldots \leq \tau_{2^{2n}}$, and then take $g_1(n) = \sigma_{2^{2n}}$ and $g_2(n) = \tau_{2^{2n}}$. In order to do this, fix an enumeration of all pairs of strings of length $n$.

Suppose that we have just defined $\sigma_i$ and $\tau_i$ and that we now need to define $\sigma_{i + 1}$ and $\tau_{i + 1}$. Let $(\gamma_1, \gamma_2)$ be the $(i + 1)^\text{th}$ pair in the enumeration fixed above. Let $\sigma_{i +1}$ and $\tau_{i+1}$ be strings extending $\sigma_i$ and $\tau_i$, respectively, such that $(\gamma_1\concat\sigma_{i + 1}, \gamma_2\concat\tau_{i + 1})$ meet all of $D_0, \ldots, D_n$ (which is possible since $D_0, \ldots, D_n$ are all dense). In order to make the choice of $\sigma_{i + 1}$ and $\tau_{i + 1}$ deterministic, choose them to first, minimize the total length $|\sigma_{i + 1}| + |\tau_{i + 1}|$ and second, be lexicographically least subject to having minimal length.

Finally, note that the entire process just described is uniformly computable from $r'$ and hence the functions $g_1$ and $g_2$ are $r'$-computable.
\end{proof}

Now inductively define a sequence of numbers $k_0, k_1, k_2, \ldots$ by setting
\begin{align*}
  k_0 &\vcentcolon= 0\\
  k_{n + 1} &\vcentcolon= k_n + |g_1(k_n)|.
\end{align*}
Note that since the function $g_1$ is $r'$-computable, so is the sequence $k_0, k_1, k_2, \ldots$. Also note that, by slightly modifying the function $g_1$ if necessary, we may assume the sequence $k_0, k_1, k_2, \ldots$ is strictly increasing.

We are now ready to construct the reals $x$ and $y$. The details of the construction are contained in the following lemma, which also captures the key property of $x$ and $y$ that we will use below.

\begin{lemma}
\label{lemma:construct}
There are reals $x, y \in \Cantor$ which are $r'$-generic such that for every $n$, either $x\restriction {[k_n, k_{n + 1})} = g_1(k_n)$ or $y\restriction {[k_n, k_{n + 1})} = g_1(k_n)$.
\end{lemma}

\begin{proof}
The idea is just that $x$ and $y$ will trade off between meeting dense sets and copying what the function $g_1$ tells them to do. In order to make this more precise, we first fix an enumeration $E_0, E_1, E_2, \ldots$ of the dense subsets of $2^{< \N}$ which need to be met to ensure $r'$-genericity.

We will define $x$ and $y$ via initial segment approximations over the course of an infinite sequence of stages. We now describe what happens on stage $n + 1$. Suppose that at the end of stage $n$ we have defined initial segments $\sigma$ and $\tau$ of $x$ and $y$, respectively. We will assume that $|\sigma| = |\tau| = k_m$ for some $m$ (and we will thus need to make sure that this holds at the end of stage $n + 1$ as well).

We proceed as follows. First, let $\sigma'$ be some string extending $\sigma$ such that $\sigma'$ meets $E_{n + 1}$. Let $m'$ be the least number such that $k_{m'} \geq |\sigma'|$. By padding $\sigma'$ with $0$s, we may assume $|\sigma'| = k_{m'}$. Let $\tau'$ be the string obtained by copying $g_1$ between $k_m$ and $k_{m'}$. In other words,
\[
  \tau' = \tau\concat g_1(k_m) \concat g_1(k_{m + 1}) \concat \cdots \concat g_1(k_{m' - 1}).
\]
Now reverse the roles of $\sigma'$ and $\tau'$ and repeat. In other words, let $\tau''$ be some string extending $\tau'$ which meets $E_{n + 1}$, let $m''$ be least such that $k_{m''} \geq |\tau''|$ and extend $\sigma'$ to $\sigma''$ by copying $g_1$ between $k_{m'}$ and $k_{m''}$. Also pad $\tau''$ with $0$s so that its length is exactly $k_{m''}$. The strings $\sigma''$ and $\tau''$ are the initial segments of $x$ and $y$, respectively, determined by the end of stage $n + 1$.
\end{proof}

We will now prove that the reals $x$ and $y$ constructed in the lemma above have the desired property. To that end, suppose $z$ is $r'$-generic. We will show that for each $n$, either the pair $(x, z)$ meets all of $D_0, D_1, \ldots, D_n$ or the pair $(y, z)$ meets all of $D_0, D_1, \ldots, D_n$. This is enough to imply that either $x$ and $z$ are mutually $r$-generic or $y$ and $z$ are mutually $r$-generic.

Fix $n$ and consider the set $F_n$ of strings of the form $\tau\concat g_2(k_m)$ where $k_m > n$ and $|\tau| = k_m$. Since $\lim_{m \to \infty}k_m = \infty$, $F_n$ is dense and since the function $g_2$ and the sequence $k_0, k_1, k_2, \ldots$ are both computable from $r'$, $F_n$ is computable from $r'$ as well. Hence any real which is $r'$-generic must meet $F_n$.

In particular, since $z$ is $r'$-generic, $z$ must meet $F_n$. Therefore, $z$ has an initial segment of the form $\tau\concat g_2(k_m)$ with $k_m > n$ and $|\tau| = k_m$. By our choice of $x$ and $y$, (see Lemma~\ref{lemma:construct} above), at least one of $x$ and $y$ has an initial segment of the form $\sigma\concat g_1(k_m)$ where $|\sigma| = k_m$. Without loss of generality, assume that $x$ has an initial segment of this form. By definition of $g_1$ and $g_2$, the pair $(\sigma\concat g_1(k_m), \tau\concat g_2(k_m))$ meets all of the dense sets $D_0, \ldots, D_n$ and hence so does the pair $(x, z)$.


\bibliographystyle{alpha}
\bibliography{bibliography}

\end{document}